\documentclass{article}
 \usepackage{amsmath}
 \usepackage{pifont}

\newtheorem{thm}{Theorem}[section]
\newtheorem{cor}[thm]{Corollary}
\newtheorem{prop}[thm]{Proposition}

\newtheorem{lem}[thm]{Lemma}
\newtheorem{Def}[thm]{Definition}

\newtheorem{ex}[thm]{Example}

\newcommand{\be}{\begin{equation}}
\newcommand{\ee}{\end{equation}}
\newcommand{\ben}{\begin{enumerate}}
\newcommand{\een}{\end{enumerate}}
\newcommand{\beq}{\begin{eqnarray}}
\newcommand{\eeq}{\end{eqnarray}}
\newcommand{\beqn}{\begin{eqnarray*}}
\newcommand{\eeqn}{\end{eqnarray*}}

\newcommand{\pa}{\partial}

\newcommand{\qed}{\hspace*{\fill}Q.E.D.}  

\begin{document}
\title{Sprays on Hamel-Funk Functions Model}
\author{Guojun Yang }
\date{}
\maketitle
\begin{abstract}
 Hamel functions of a spray play an important role in the study of the projective metrizability of the concerned spray,
 and Funk functions are special Hamel functions. A
 Finsler metric is a special Hamel function of the spray induced by the
 metric itself and a Funk metric is a special Funk function of a
 Minkowski spray. In this paper, we study sprays on a Hamel or Funk function model.
 Firstly, we give some basic properties of a Hamel or Funk
 function of a spray and some curvature properties of  a Hamel or Funk
 function in projective relations. We use the Funk metric to construct
 a family of sprays and obtain some of their curvature properties and their metrizability conditions.
  Secondly, we consider the
 existence of Funk functions on certain spray manifold. We prove
 that there exist local Funk functions on a R-flat spray manifold,
 and on certain projectively flat Berwald spray manifolds, we construct a multitude of nonzero Funk
 functions. Finally, we introduce a new class of sprays called  Hamel or Funk sprays
  associated to given sprays and  Hamel or Funk functions. We obtain some
 special properties of a Hamel or Funk spray of scalar curvature,
 especially on its  metrizability and a special form
 of its Riemann curvature.

 \


\end{abstract}

\section{Introduction}

Spray geometry  studies the properties of  sprays on a manifold,
and it is closely related to Finsler geometry. Every Finsler
metric induces a natural spray but there are a lot of sprays which
cannot be induced by any Finsler metric (\cite{BM2, EM, LMY,
Yang1, Yang2}). A spray ${\bf G}$ on a manifold $M$ defines a
special vector filed on a conical region $\mathcal{C}$ of
 $TM\setminus \{0\}$, and  it naturally defines the integral curves and the
 projections of the integral curves onto the manifold $M$ are
 called geodesics. For deep understanding of spray geometry, it is
 important to develop more basic curvatures, to study some classes of sprays with
 special geometric properties and to investigate the (projective) metrizablity
 of a spray under certain curvature conditions.

Many basic curvatures, such as Riemann curvature, Ricci curvature,
Weyl curvature, Berwald curvature and Douglas curvature, appearing
in Finsler geometry, are actually defined in spray geometry via
the spray coefficients. With a fixed volume form, the S-curvature,
originally introduced in Finsler geometry by Z. Shen in
\cite{Shen0}, can also be generalized to spray geometry. The
notions for a Finsler metric of scalar flag curvature, isotropic
flag curvature and constant flag curvature, already have the
corresponding generalization in spray geometry, and they are
respectively called a spray of scalar curvature, isotropic
curvature and constant  curvature. A spray of scalar curvature
first appears in Shen's book \cite{Shen6} (therein, it is
originally called an isotropic spray). In \cite{LS}, B. Li and Z.
Shen introduce a spray of isotropic curvature. The present author
introduces a spray of constant curvature in \cite{Yang2}. Further
studies show that these curvatures or curvature characters play an
important role in spray geometry.

Among all sprays, some special classes of sprays are helpful in
understanding the geometric structures of a general spray
(\cite{BM2, LMY, Yang2}). A Berald spray means that its spray
coefficients are polynomials and every Berwald metric induces a
Berwald spray. A (locally) projectively flat spray means its
geodesics are (local) straight lines, and it's the Hilbert's
Fourth Problem to study the structure  of a (locally) projectively
flat spray, especially when the spray is induced by a Finsler
metric. These sprays and sprays of scalar curvature (including
isotropic, constant curvature) have some special good properties
and deserve further studies.

The metrizability problem for a spray ${\bf G}$ seeks for a
Finsler metric whose spray is just {\bf G}, and the projective
metrizability of a  spray {\bf G} aims to look for a Finsler
metric projectively related to {\bf G}. So a natural question is
to determine whether a given spray is (projectively)
Finsler-metrizable or not under certain curvature conditions.
Quite a few papers concentrate on this problem (\cite{BM,  BM2,
Cram, EM, LMY, Ma1,  Mu, Yang1, Yang2}). For instance, any spray
of scalar curvature is locally projectively metrizable (\cite{Ma1,
BM}); the metrizability for a spray of scalar curvature with
nonzero Ricci curvature  is solved (\cite{BM2});
 the local structure is determined
for a projectively flat Berwald spray of isotropic curvature to be
metrizable (\cite{Yang2}).

In this paper, we are going to study some geometric properties of
Hamel-Funk functions and Hamel-Funk sprays. The notion of Funk
function first appears in Shen's book \cite{Shen6}.  Hamel-Funk
sprays are a new class of sprays  we introduce in this paper.

Let {\bf G} be a spray on a manifold $M$. We can use the Berwald
connection of {\bf G} to define a Hamel function or a Funk
function of {\bf G} (Definition \ref{DHF1}). When a Finsler metric
$F$ is projectively flat, $F$ is a Hamel function of the spray
${\bf G}=0$, which satisfies the well-known Euler-Lagrange
equations. When $F$ is a Funk metric on a convex domain in $R^n$,
$F$ is a Funk function of the spray ${\bf G}=0$. In Sec.
\ref{sec3}, we introduce some basic properties of a Hamel or Funk
function of a Minkowski spray or a general spray, and generalize
some known equations  in \cite{Shen5} for projectively flat
Finsler metrics of constant flag curvature to Hamel or Funk
functions of a spray (Lemmas \ref{Plem27}--\ref{Plem30}). We see
that a Funk function is a Hamel function and the collection of all
Hamel functions forms a linear space (Lemma \ref{Plem25}).

 A Hamel or
Funk function plays a special role in projective relations. A
simple fact shows that a Hamel function is projectively invariant,
but a Funk function is not (Lemmas \ref{PPHF1}, \ref{PPHF01}). Z.
Shen proves that The Riemann curvatures of two projectively
related sprays keep invariant if the projective factor is a Funk
function (\cite{Shen6}). Further, we prove that if one of two
projectively related sprays is of isotropic curvature, then the
other is also of isotropic curvature if and only if the projective
factor if a Hamel function (Proposition \ref{PFH1}(iii)). Starting
from a Funk metric, we obtain the following result.

\begin{thm}\label{MTh1}
 Let {\bf G} be the spray of the Funk metric $F$ on a convex domain $\Omega\subset R^n$.  For a
constant $c$, define a spray
 \be\label{ME1}
 \bar{G}^i(y):=G^i(y)+cF(-y)y^i.
 \ee
 \ben
  \item[{\rm (i)}] ${\bf \bar{G}}$ is  of isotropic curvature.
  Further, ${\bf \bar{G}}$  is of constant curvature if and only
  if $c=0$ or $c=-1/2$.
 \item[{\rm (ii)}] ${\bf \bar{G}}$ is (locally) metrizable if and
only if $c=0$ or $c=-1/2$. In this case,  ${\bf \bar{G}}$ is
induced by $F$ if $c=0$; and if $c=-1/2$, ${\bf \bar{G}}$ is
induced by the Klein metric
 \be\label{Klein}
\bar{F}(y):=\big[F(y)+F(-y)\big]/2.
 \ee
 \een
\end{thm}

Theorem \ref{MTh1}(i)  gives a new class of examples to show that
a spray of isotropic curvature is not necessarily of constant
curvature even in dimension greater than two (cf. \cite{Yang2}).
Theorem \ref{MTh1}(ii) also gives a class of non-metrizable
sprays. For ${\bf \bar{G}}$ in Theorem \ref{MTh1}, if $F(-y)$ is
replaced by $F(y)$, the result refers to \cite{Yang1}.

By some properties of a Hamel function, it is relatively free to
construct  Hamel functions of certain sprays. For instance, a
Finsler metric is a Hamle function of its induced spray, and
Propositions \ref{PEHF001} and \ref{PEHF002} also give some Hamel
functions of certain sprays. However, the construction of a Funk
function of a  spray is not easy. On a Minkowski spray ${\bf
G}=0$, we can get many Funk functions (a Funk metric itself is).
In \cite{Shen6}, Z. Shen proves that no nonzero Funk function
exists on a compact regular spray manifold (Theorem \ref{TEHF1}).
In \cite{BM0}, I. Bucataru and Z. Muzsnay show that no nonzero
Funk function exists on a Finsler manifold of nonzero scalar
curvature (Theorem \ref{TEHF2}). To construct nontrivial Funk
functions of a spray, we first prove Proposition \ref{PEHF1} on
some properties satisfied by a Funk function of a spray of scalar
or isotropic or constant curvature, and then we give many
nontrivial closed 1-form Funk functions of projectively flat
Berwald sprays in Examples \ref{EEHF1}--\ref{EEHF3}. Further we
have the following local existence theorem.

\begin{thm}\label{MTh2}
 There exist non-trivial (local) Funk functions on a R-flat spray
 manifold.
\end{thm}

A locally Minkowski spray is R-flat and there are many Funk
functions of a Minkowski spray (see  Sec. \ref{sec3}). If a spray
is just R-flat (not Berwaldian), we generally don't have the
explicit expression of the Funk function in Theorem \ref{MTh2}.

Using a  spray {\bf G} and a nonzero homogeneous function $Q$ of
degree one, we define a spray ${\bf \bar{G}}$ by
 \be\label{YHFS1}
\bar{G}^i:=G^i+\frac{Q_{;0}}{2Q}y^i.
 \ee
Every spray ${\bf \bar{G}}$ in the projective class of {\bf G} can
be locally written in the form (\ref{YHFS1}) for some $Q$ (Lemma
\ref{LHFS001}), and it is projectively invariant for a fixed $Q$
(Lemma \ref{LHFS01}). When $Q$ is a Hamel or Funk function, we
obtain  a Hamel or Funk spray ${\bf \bar{G}}$  associated to
$({\bf G},Q)$, with ${\bf G},Q$ being the adjoint spray and the
adjoint Hamel function respectively (Definition \ref{DHFS63}).
Hamel functions and Hamel sprays play an important role in
studying the projective metrizability of a given spray. We have
the following theorem.

\begin{thm}\label{THF1}
  A spray {\bf G} is
projectively metrizable if and only if  it has a Finsler-Hamel
function $\bar{F}$. In this case, the spray {\bf G}  is
projectively induced by the Finsler metric $\bar{F}$ whose spray
is the Hamel spray associated to $({\bf G},\bar{F})$.
\end{thm}

For the projective metrizability of a spray, a different
description from Theorem \ref{THF1} is given in \cite{BM}. In
Theorem \ref{THF1}, if the Hamel function is a Finsler metric,
then the Hamel spray is a Finsler spray. As we know, for a Finsler
spray of scalar curvature, its Riemann curvature has a special
form (\ref{SCFF}). This form is inherited well by the Riemann
curvature of a Hamel spray of scalar curvature in a similar way.

\begin{thm}\label{MTh4}
{\rm(i)} A Hamel spray ${\bf \bar{G}}$  is of scalar curvature if
and only if its Riemmann curvature has the form
 \be\label{YHFS4}
 \bar{R}^i_{\ k} =\lambda (Q^2\delta^i_k-QQ_{.k}y^i),
 \ee
 where $\lambda=\lambda(x,y)$ is a scalar function and $Q$ is the
 adjoint Hamel function of ${\bf \bar{G}}$.

 {\rm (ii)} A Funk spray ${\bf \bar{G}}$  is of isotropic curvature if
and only if its Riemmann curvature has the form (\ref{YHFS4}),
where $\lambda=\lambda(x)$ is a scalar function and $Q$ is the
 adjoint Funk function of ${\bf \bar{G}}$.
\end{thm}

The scalar function $\lambda$ in Theorem \ref{MTh4} is determined
by (\ref{HFS01}) and (\ref{HFS02}) respectively, and $\lambda$ can
be a constant for some special cases (cf. Examples \ref{EHFS2} and
\ref{EHFS3}).

A given Hamel spray is possibly associated to $({\bf G},Q_1)$ and
$({\bf G},Q_2)$ with different Hamel functions $Q_1,Q_2$ of {\bf
G}. But under certain condition, a given Hamel spray corresponds a
unique Hamel function (up to a constant scaling) (cf. Proposition
\ref{PHFS3}). On  the metrizability of a Hamel spray, we have the
following theorem.

\begin{thm}\label{THFS02}
 A Hamel spray of nonzero scalar curvature is metrizable if and only if its adjoint Hamel
 function is a Finsler metric.  In this case, the Hamel spray is
 uniquely
 induced by its adjoint Hamel
 function (the metric is unique up to a constant scaling).
\end{thm}

 Theorem \ref{THFS02} is not true if the Hamel
spray is of zero scalar curvature or not of scalar curvature. A
Minkowski spray ${\bf \bar{G}}$ is a natural Hamel spray of zero
scalar curvature, but any two Minkowski metrics induce ${\bf
\bar{G}}$. Any Finsler spray ${\bf \bar{G}}$ (induced by a Finsler
metric $\bar{F}$) is also a natural Hamel spray associated to
$({\bf \bar{G}},\bar{F})$. If ${\bf \bar{G}}$ is not of scalar
curvature, we can find many Finsler metrics inducing the same
${\bf \bar{G}}$ but not having constant ratio (see \cite{CS}:
Chapter 9).

Starting from a spray {\bf G} of scalar curvature and  a closed
1-form, we can obtain many Hamel sprays which are not metrizable
(cf. Examples \ref{EHFS2} and \ref{EHFS3}, Corollary
\ref{CHFS001}).

\section{Preliminaries}\label{pre}

Let $M$ be an $n$-dimensional manifold. A conical region
$\mathcal{C}=\mathcal{C}(M)$ of $TM\setminus \{0\}$ means
$\mathcal{C}_x:=\mathcal{C}\cap T_xM\setminus \{0\}$ is conical
region for $x\in M$ ($\lambda y\in \mathcal{C}_x$ if
$\lambda>0,y\in \mathcal{C}_x$). A {\it spray} on $M$ is a smooth
vector field ${\bf G}$ on a conical region $\mathcal{C}$ of
 $TM\setminus \{0\}$ expressed in
a local coordinate system $(x^i,y^i)$ in $TM$ as follows
 $${\bf G}=y^i\frac{\pa}{\pa x^i}-2G^i\frac{\pa}{\pa y^i},$$
 where $G^i=G^i(x,y)$ are local functions which are positively homogeneous with degree
 two. If {\bf G} is defined on the whole $TM\setminus \{0\}$, we
 call {\bf G} a regular spray.

The  Riemann curvature tensor $R^i_{\ k}$ of a given spray $G^i$
is defined by
 \be\label{y004}
  R^i_{\ k}:=2\pa_k G^i-y^j(\pa_j\dot{\pa_k}G^i)+2G^j(\dot{\pa_j}\dot{\pa_k}G^i)-(\dot{\pa_j}G^i)(\dot{\pa_k}G^j),
 \ee
where we define $\pa_k:=\pa/\pa x^k,\dot{\pa}_k:=\pa/\pa y^k$. The
trace of $R^i_{\ k}$ is called the Ricci curvature, $ Ric:=R^i_{\
i}$. A spray ${\bf G}$ is said to be {\it R-flat} if $R^i_{\
k}=0$. A spray ${\bf G}$ is said to be of {\it scalar curvature}
if  its Riemann curvature $R^i_{\ k}$ satisfies
 \be\label{sc}
  R^i_{\ k}=R\delta^i_k-\tau_ky^i,
  \ee
where $R=R(x,y)$ and $\tau_k=\tau_k(x,y)$ are some homogeneous
functions (\cite{Shen6}). If in (\ref{sc}) there holds
$R_{.i}=2\tau_i$,
 then ${\bf G}$  is said to be of {\it isotropic curvature}
 (\cite{LS}). A spray ${\bf G}$ is said to be of
 {\it constant curvature} if ${\bf G}$ satisfies (\ref{sc}) with
 (\cite{Yang2})
  $$
 \tau_{i;k}=0 \ (\ \Leftrightarrow \ R=\tau_k=0,\ or \ R_{;i}=0).
  $$
The three notions are natural generalizations for a Finsler metric
of scalar/isotropic/constant flag curvature in Finsler geometry.

 A spray {\bf G} is called a Berwald spray if its Berwald curvature vanishes
 $G^i_{hjk}:=\dot{\pa}_h\dot{\pa}_j\dot{\pa}_kG^i=0$.  If $G^i$ is locally given by $G^i=0$, we call {\bf
G} a locally Minkowski spray. Two sprays ${\bf G}$ and $\bar{\bf
G}$ on $M$ are said to be  projectively
 related if their geodesics coincide as  set points, or
 equivalently,  $\bar{G}^i=G^i+Py^i$
 for some positively homogeneous function $P$ of degree one.
 A spray {\bf G} is said to be locally projectively flat if
 locally $G^i$ can be expressed as
$G^i=Py^i$. A spray {\bf G} is called {\it reversible}  if the
spray $G^i(-y)$ is projectively related to $G^i(y)$, and {\it
strictly reversible} if $G^i(y)=G^i(-y)$. We
 use $Proj({\it G})$ to denote the projective class of {\bf G}.

In the calculation of some geometric quantities of a spray, it is
very convenient to use Berwald connection as a tool. For a spray
manifold $({\bf G}, M)$, Berwald connection is usually defined as
a linear connection on the pull-back $\pi^*TM$ ($\pi:
TM\rightarrow M$ the natural projection) over the base manifold
$M$. The Berwald connection is defined by
  $$
  D(\pa_i)=(G^k_{ir}dx^r)\pa_k,\ \ \ \ \ (G^k_{ir}:=\dot{\pa}_r\dot{\pa}_iG^k),
  $$
 For a spray tensor $T=T_idx^i$ as an example,   the horizontal and vertical derivatives of $T$ with respect to Berwald
 connection are given by
  $$
 T_{i;j}=\delta_jT_i-T_rG^r_{ij},\ \ \ \ \ \ \
 T_{i.j}=\dot{\pa}_jT_i,\ \ \ \ (\delta_i:=\pa_i-G^r_i\dot{\pa}_r).
  $$
  For the Ricci identities and Bianchi identities of Berwald connection, readers
 can refer to \cite{AIM}.

 In this paper, we  define a Finsler metric $L(\ne 0)$ on a manifold $M$ as follows (cf. \cite{Shen6}):
 (i) for any $x\in M$, $L_x$ is defined on a conical region of  $T_xM\setminus \{0\}$
  and $L$ is $C^{\infty}$; (ii) $L$
 is positively homogeneous of degree two; (iii) the fundamental
metric tensor $g_{ij}:=(\frac{1}{2}L)_{y^iy^j}$ is non-degenerate.
 A Finsler metric
$L$ is said to be regular if additionally $L$ is defined on the
whole  $TM\setminus \{0\}$ and $(g_{ij})$ is positively definite.
Otherwise,
 $L$ is called singular. In general case, we
don't require that $L$ be regular. If a Finsler metric $L>0$, we
put $L=F^2$, and in this case, $F$ is also called a Finsler metric
being positively homogeneous of degree one.

 Any Finsler metric $L$ induces
a natural spray whose coefficients $G^i$ are given by
 $$
 G^i:=\frac{1}{4}g^{il}\big \{L_{x^ky^l}y^k-L_{x^l}\big
 \},
$$
where $(g^{ij})$ is the inverse of $(g_{ij})$. $L$ is said to be
of {\it scalar flag curvature} $K=K(x,y)$ if
 \be\label{SCFF}
 R^i_{\ k}=K(L\delta^i_k-y^iy_k)=K(F^2\delta^i_k-FF_{.k}y^i), \ \
 \ \ (L=F^2),
 \ee
 where $y_k:=(L/2)_{.k}=g_{km}y^m$. If $K_{.i}=0$, then $L$ is said to be of
 {\it isotropic flag curvature}, and in this case, $K$ is a
 constant if the dimension $n\ge 3$ (Schur's Theorem). $L$ is said to be of
 {\it constant flag curvature} if $K$ is constant.

  A spray {\bf G} is (globally)
 Finsler-metrizable on $M$ (or on $\mathcal{C}(M)$) if there is a Finlser metric $L$ defined on
 a conical region $\mathcal{C}(M)$ and $L$ induces {\bf G}. A spray {\bf G} is
 locally Finsler-metrizable on $M$ if for each $x\in M$, there is
 a neighborhood $U$ of $x$ such that {\bf G} is
 Finsler-metrizable on $U$.

\begin{lem}\label{Prelem1}
 A Finsler spray
 ${\bf \bar{G}}$ and a spray {\bf G} are related by
  \be\label{pr0}
 \bar{G}^i=G^i+\frac{\bar{F}_{;0}}{2\bar{F}}y^i-\frac{1}{2}\bar{F}\bar{g}^{ir}(\bar{F}_{;r}-\bar{F}_{.r;0}),
  \ee
  where $\bar{F}$ is a Finsler metric inducing ${\bf \bar{G}}$.
\end{lem}

\begin{lem}\label{Prelem2}
 A Finsler spray ${\bf \bar{G}}$ (induced by $\bar{F}$) is
 projective to a spray {\bf G} if and only if $\bar{F}$ satisfies
 one of the following conditions:
  \be\label{pr3}
 {\rm (i)}\ \bar{F}_{;i}=(P\bar{F})_{.i},\ \ \ \ \  {\rm (ii)}\  \bar{F}_{;i}=\bar{F}_{.i;0},
 \ \ \ \ \  {\rm (iii)} \ \bar{F}_{;i}=\frac{1}{2}\bar{F}_{;0.i},
  \ee
where in case (i), $P$ is a homogeneous function of degree one. In
this case, $P$ is  the projective factor which is given by
   \be\label{pr4}
 P=\frac{\bar{F}_{;0}}{2\bar{F}}.
   \ee
\end{lem}

\begin{lem}
 Let ${\bf \bar{G}}$ and {\bf G} be projectively related with
 $\bar{G}^i=G^i+Py^i$. Then their Riemann curvatures and Ricci
 curvatures are related by
  \beq
 \bar{R}^i_{\ k}&&\hspace{-0.6cm}=R^i_{\ k}+(3Q_k-Q_{0.k})y^i-Q_0\delta^i_k,\label{pr9}\\
 \bar{R}ic&&\hspace{-0.6cm}=Ric-(n-1)Q_0=Ric-(n-1)(P_{;0}-P^2),\label{pr13}\\
 &&\hspace{-1.8cm}\big(Q_j:=P_{;j}-PP_j,\ \ \ \
 Q_{jk}:=P_{j;k}-P_{k;j}=Q_{k.j}-Q_{j.k}\big).\nonumber
  \eeq
\end{lem}

The above  three lemmas can be found  in some textbooks. For the
proofs, readers can refer to Chapter 3 of \cite{AIM}, for
instance.

\begin{lem}(\cite{Yang2})\label{Yplem}
Let $T$ be a positively homogeneous function of degree zero on a
spray manifold of scalar curvature with nonzero Ricci curvature.
If $T_{;i}=0$, then $T$ is  constant.
\end{lem}

\begin{lem}\label{plem1}
 On a pray manifold $(M,{\bf G}$), put
  \be
  T_i:=Q_{;r.i}y^r-Q_{;i}\ (=Q_{.i;0}-Q_{;i}),\label{p018}
  \ee
  where $Q$ is an arbitrary positively homogeneous function of degree
  one. Then $T_i$ is a projective invariant.
\end{lem}

{\it Proof :} For two projectively related ${\bf \bar{G}}$ and
{\bf G}, let $\bar{G}^i=G^i+Py^i$. A direct computation gives
 $$
 Q_{\bar{;}i}=Q_{;i}-(PQ)_{.i},\ \ \ \ \
 Q_{\bar{;}r.i}=Q_{;r.i}-(PQ)_{.i.r}.
 $$
Therefore, it is easy to get
$$Q_{\bar{;}r.i}y^r-Q_{\bar{;}i}=Q_{;r.i}y^r-Q_{;i},$$
which shows that $T_i$ is projectively invariant.   \qed

\section{Hamel-Funk Functions and Basic Properties}\label{sec3}

In this section, we introduce the definitions of a Hamel function
and a Funk function, and then study some of their basic properties
on  a  spray manifold.

\begin{Def}\label{DHF1}
 Let $(M,{\bf G}$) be a spray manifold. A  positively homogeneous function $P$  of degree
  one  is called a Hamel
  Function (resp. Funk function, weak Funk function) of {\bf G} if $P$ satisfies
   $$
 P_{;k}=P_{k;0} \ (resp. \ P_{;k}=PP_k,\  P_{;0}=P^2),\ \ \ (P_k:=P_{.k}).
   $$
\end{Def}

Z. Shen introduces Funk function and weak Funk function of a spray
in \cite{Shen6}.

\subsection{Hamel-Funk Functions of a Minkowski Spray}
In this subsection,  we show some properties of a Hamel-Funk
function of the Minkowski spray ${\bf G}=0$ on $R^n$ (cf.
\cite{CS, LS, Shen5}).

\begin{lem}\label{Plem22}
 Let $P$ be a Funk function of the spray ${\bf G}=0$ on
 $R^n$. Then we have
  \beq
&& P_{x^{i_1}\cdots
 x^{i_m}}=\frac{1}{m+1}\big(P^{m+1}\big)_{y^{i_1}\cdots
 y^{i_m}},\label{PHF01}\\
&&Q_{x^k}=(PQ)_{y^k}\ \ \ \ \Longrightarrow \ \ \ \
Q_{x^{i_1}\cdots x^{i_m}}=(QP^m)_{y^{i_1}\cdots
 y^{i_m}},\label{PHF02}
   \eeq
   where in (\ref{PHF02}), $Q$ is some function on $TR^n$.
\end{lem}

\begin{lem}\label{HF2}
 Let $P$ be a Funk function of the spray ${\bf G}=0$  on $R^n$ and $k\ge 0$ be an integer. Define
  $$
 Q:=(P^{k+1})_{y^{i_1}\cdots y^{i_k}}x^{i_1}\cdots x^{i_k}\ \ \Big(=(k+1)P_{x^{i_1}\cdots
 x^{i_k}}x^{i_1}\cdots x^{i_k}\Big).
  $$
 Then $Q$ is a Hamel
 function of {\bf G}.
\end{lem}

{\it Proof :} This proof refers to \cite{LS}. By
$P_{x^k}=PP_{y^k}$, we have
 \beqn
 Q_{x^j}&&\hspace{-0.6cm}=(k+1)(P^kP_{x^j})_{y^{i_1}\cdots y^{i_k}}x^{i_1}\cdots
 x^{i_k}+k(P^{k+1})_{y^{i_1}\cdots y^{i_{k-1}}y^j}x^{i_1}\cdots
 x^{i_{k-1}}\\
  &&\hspace{-0.6cm}=\frac{k+1}{k+2}(P^{k+2})_{y^{i_1}\cdots
  y^{i_k}y^j}x^{i_1}\cdots x^{i_k}+k(P^{k+1})_{y^{i_1}\cdots y^{i_{k-1}}y^j}x^{i_1}\cdots
 x^{i_{k-1}}.
 \eeqn
 Then we easily obtain
  $$Q_{x^j}=Q_{x^my^j}y^m,$$
  which shows that $Q$ is a Hamel
 function of {\bf G}.         \qed

\begin{lem}(\cite{CS})\label{Plem24}
 Let $\varphi=\varphi(y)$ be a positively homogeneous function of
 degree one on
 $R^n$ and be $C^{\infty}$ on $R^n\setminus \{0\}$. Then there is a neighborhood $U$ of\ \ $0\in R^n$, and a unique Funk function
 $P=P(x,y)$ of the spray ${\bf G}=0$ such that for $x\in U$, $P$ satisfies
  \be\label{prof1}
 P(x,y)=\varphi(y+P(x,y)x).
  \ee
  In this case, we have $P(0,y)=\varphi(y)$.
\end{lem}

By Lemma \ref{Plem24}, we can construct many Funk functions of the
spray ${\bf G}=0$ by solving $P$ from (\ref{prof1}) for a given
$\varphi$. For example, let $\varphi(y)=|y|$ and
$\varphi(y)=-\langle a,y\rangle$ respectively, where $a$ is a
constant vector. Then solving $P$ from (\ref{prof1}) gives
respectively
 \be
 P=\frac{\sqrt{(1-|x|^2)|y|^2+\langle x,y\rangle^2}+\langle
 x,y\rangle}{1-|x|^2},\ \ \ \ \ P=-\frac{\langle a,y\rangle}{1+\langle
 a,x\rangle}.
 \ee

\subsection{Hamel-Funk Functions of a General Spray}

\begin{lem}\label{Plem25}
 A Hamel function or a Funk function has the following basic
properties:
 \ben
  \item[{\rm (i)}] For a spray, any Funk function is a Hamel function.
  \item[{\rm (ii)}] The collection of Hamel
functions on  a spray manifold is a vector space under the natural
addition and number-product for functions. So any linear
combination of any two Funk functions is a Hamel function.
 \een
\end{lem}

{\it Proof :}  For case (i), if $P$ is a Funk function of a spray,
then we have $P_{;k}=PP_k$, which gives $P_{;0}=P^2$. So we have
$P_{;0.k}=2PP_k$, that is, $P_{;k}+P_{k;0}=2PP_k$. Thus we obtain
$P_{;k}=P_{k;0}\ (=PP_k)$. This means that $P$ is a Hamel
function. For case (ii), the former is obvious by the definition
of a Hamel function and the latter follows from case (i).
 \qed

\begin{lem}\label{Plem26}
 Let {\bf G} be a spray and $Q$ is a Hamel function of {\bf G}. Put
 $P(y):=Q(-y)$.
 Then $P$ is a Hamel function of {\bf G} if and only if
  \be\label{cwy75}
 \big[G^r(y)-G^r(-y)\big]Q_{.r.i}(-y)=0.
  \ee
  In particular, (\ref{cwy75}) holds if {\bf G} is locally
  projectively flat or  reversible.
\end{lem}

{\it Proof :} We first have
 \beq
 P_{;i}(y)\hspace{-0.6cm}&&=(\pa_iP-G^r_iP_{.r})(y)=\pa_iQ(-y)+G^r_i(y)Q_{.r}(-y)\nonumber\\
\hspace{-0.6cm}&&=Q_{;i}(-y)+\big[G^r_i(y)+G^r_i(-y)\big]Q_{.r}(-y).\label{cwy77}
 \eeq
 Next we get
 $$
 P_{;m.i}(y)=-Q_{.i;m}(-y)+\big[G^r_{mi}(y)-G^r_{mi}(-y)\big]Q_{.r}(-y)
 -\big[G^r_m(y)+G^r_m(-y)\big]Q_{.r.i}(-y).
 $$
Therefore, we obtain
 $$
 P_{i;0}(y)=Q_{.i;0}(-y)+\big[G^r_{i}(y)+G^r_{i}(-y)\big]Q_{.r}(-y)-2\big[G^r(y)-G^r(-y)\big]Q_{.r.i}(-y).
 $$
Since $Q$ is a Hamel function, we have
 $$
(P_{;i}-P_{i;0})(y)=2\big[G^r(y)-G^r(-y)\big]Q_{.r.i}(-y),
 $$
which shows that $P$ is a Hamel function of {\bf G} if
(\ref{cwy75}) holds. \qed

\

Lemma \ref{Plem26} shows that  the Klein metric in (\ref{Klein})
is a Hamel function and it is projectively flat.

In \cite{Shen5}, Z. Shen gives some equations satisfied by a
projectively flat Finsler metric of constant flag curvature and
its corresponding projective factor. These properties can be
generalized to a general spray space satisfied by a Hamel or Funk
function (see Lemmas \ref{Plem27}--\ref{Plem30}).

\begin{lem}\label{Plem27}
 Let $Q$ be a Funk function of a spray {\bf G} and put
 $ \widetilde{Q}(y):=\big[Q(y)+Q(-y)\big]/2.$
 Then $Q$ satisfies
  \be\label{P315}
 P:=\frac{Q_{;0}}{2Q}=\frac{1}{2}Q,\ \ \ \ \
 \frac{P^2-P_{;0}}{Q^2}=-\frac{1}{4}.
  \ee
If {\bf G} is strictly reversible, then $\widetilde{Q}$ is a Hamel
function of {\bf G} and  satisfies
 \be\label{P316}
 \widetilde{P}:=\frac{\widetilde{Q}_{;0}}{2\widetilde{Q}}=\frac{1}{2}\big[Q(y)-Q(-y)\big],\ \ \ \ \
 \frac{\widetilde{P}^2-\widetilde{P}_{;0}}{\widetilde{Q}^2}=-1.
  \ee
\end{lem}

{\it Proof :} We only prove (\ref{P316}). Since {\bf G} is
strictly reversible, we have $G^i(y)=G^i(-y)$. Then we get
$G^r_i(-y)=-G^r_i(y)$. By the calculation of (\ref{cwy77}) and the
definition of a Funk function, we have
 \beqn
 &&\widetilde{P}=\frac{1}{2}\frac{Q_{;0}(y)-Q_{;0}(-y)}{Q(y)+Q(-y)}
 =\frac{1}{2}\frac{Q^2(y)-Q^2(-y)}{Q(y)+Q(-y)}=\frac{1}{2}\big[Q(y)-Q(-y)\big],\\
 &&\widetilde{P}^2-\widetilde{P}_{;0}=\frac{1}{4}\big[Q(y)-Q(-y)\big]^2-\frac{1}{2}\big[Q^2(y)+Q^2(-y)\big]
 =-\widetilde{Q}^2.
 \eeqn
By Lemma \ref{Plem26}, $Q(-y)$ is a Hamel function of {\bf G}.
Thus $\widetilde{Q}$ is a Hamel function of {\bf G} by Lemma
\ref{Plem25}. This gives the proof. \qed

\begin{lem}\label{Plem28}
Let $Q$ and $\bar{Q}$ be two Funk functions of a spray {\bf G} and
put
 $\widetilde{Q}:=(Q-\bar{Q})/2.$
 Then $\widetilde{Q}$ is a Hamel function of {\bf G} and satisfies
\be\label{P317}
 \widetilde{P}:=\frac{\widetilde{Q}_{;0}}{2\widetilde{Q}}=\frac{1}{2}(Q+\bar{Q}),\ \ \ \ \
 \frac{\widetilde{P}^2-\widetilde{P}_{;0}}{\widetilde{Q}^2}=-1.
  \ee
\end{lem}

It is easy to verify that  a Funk function $Q$  of a spray {\bf G}
satisfies
 \be\label{PY313}
  Q_{;i.j}=Q_{;j.i}.
  \ee

\begin{lem}\label{Plem29}
Let $P,Q$ be two positively homogeneous
 functions of degree one satisfying
  \be\label{PY314}
  Q_{;i}=(PQ)_{.i}
  \ee
with respect to a spray {\bf G}.  Then $Q$
 is a Hamel function of {\bf G}. Further, if ${\bf G}=0$ and $P$ is a Funk function of ${\bf
 G}$, then in a neighborhood of $0\in R^n$, $Q$ can be
 expressed as
  \be\label{PY315}
 Q(x,y)=\psi\Big(y+P(x,y)x\Big) \big(1+P_{y^m}(x,y)x^m\big),
  \ee
  where $\psi$ is defined by $\psi(y):=Q(0,y)$.
\end{lem}

{\it Proof :} It is easy to prove that $Q$
 is a Hamel function of {\bf G} by (\ref{PY314}). Now let ${\bf G}=0$ and $P$ be a
 Funk function of ${\bf G}$ satisfying  (\ref{PY314}), that is,
 $Q_{x^k}=(PQ)_{y^k}$.
 Then by (\ref{PHF02})
 we have
  \be\label{PY316}
 Q_{x^{i_1}\cdots x^{i_m}}(0,y)=(\psi\varphi^m)_{y^{i_1}\cdots
 y^{i_m}}(y),\ \ \ \ \ (\varphi(y):=P(0,y)).
  \ee
On the other hand, if $Q$ is given by (\ref{PY315}), then using
$P_{x^k}=PP_{y^k}$ and $P_{x^ky^m}=P_{x^my^k}$ (see
(\ref{PY313})), we can directly verify $Q_{x^k}=(PQ)_{y^k}$, which
gives (\ref{PY316}) again. So by the uniqueness of expansion, $Q$
is determined by (\ref{PY315}) in a neighborhood of $0\in R^n$.
\qed

\begin{lem}\label{Plem30}
 Let $H=P+iQ$ be a complex-valued Funk  function of a spray {\bf
 G}. Then $P,Q$ are related by
  \be\label{PY317}
 Q_{;k}=(PQ)_{.k},\hspace{1cm}  P_{;k}=PP_{.k}-QQ_{.k}.
  \ee
  So $Q$ is a Hamel function.
\end{lem}

Let $\phi=\phi(y)$ be a complex-valued homogeneous function of
degree one on $R^n$. If a complex-valued function $H=H(x,y)$
satisfies
 \be\label{PY318}
 H(x,y)=\phi\big(y+H(x,y)x\big),
 \ee
then $H$ is a Funk function, $H_{x^k}=HH_{y^k}$, of ${\bf G}=0$ on
$R^n$ (see  \cite{Shen5}).

\begin{ex}
 Let $\phi(y)=i|y|$ in (\ref{PY318}), that is,
  $$
 H=i\sqrt{|y|^2+2\langle x,y\rangle H+|x|^2H^2}.
  $$
  Solving this equation yields $H=P+iQ$ satisfying (\ref{PY317})
  with
   $$
 P=-\frac{\langle x,y\rangle }{1+|x|^2},\ \ \ \
 Q=\frac{\sqrt{(1+|x|^2)|y|^2-\langle x,y\rangle^2}}{1+|x|^2}.
   $$
\end{ex}

\section{Hamel-Funk Functions in Projective Relations}\label{sec4}

Now  we investigate some projective properties related to a Hamel
or Funk function.

\begin{lem} \label{PPHF1}
  A Hamel function of a spray is
projectively invariant.
\end{lem}

Lemma \ref{PPHF1} directly follows from Lemma \ref{plem1}. A Funk
function is not projectively invariant, as shown in the following
lemma.

\begin{lem}\label{PPHF01}
 If two projectively related sprays have a same nonzero Funk
 function, then the two sprays are equal.
\end{lem}

{\it Proof :} Let ${\bf
 \bar{G}}$ and ${\bf G}$ be related by $\bar{G}^i=G^i+Py^i$. For a positively homogeneous function $Q\ne 0$ of degree one, we have
 \be\label{proj1}
Q_{\bar{;}k}=Q_{;k}-(PQ)_{.k},
 \ee
 which implies  $PQ=0$ if $Q$ is a Funk function of {\bf G} and also
 of ${\bf \bar{G}}$.    \qed

\begin{lem}\label{PPHF001}
 For ${\bar{\bf G}}\in Proj({\bf G})$ with $\bar{G}^i=G^i+Py^i$,
 $P$ is a Funk function of {\bf G} if and only if
 $P_{\bar{;}k}=-PP_k$, and $P$ is a weak Funk function of {\bf G} if and only if
 $P_{\bar{;}0}=-P^2$.
\end{lem}

Lemma \ref{PPHF001} is a direct result of (\ref{proj1}) by taking
$Q=P$. Now we consider some curvature properties  when the
projective factor is a Hamel of Funk function.

\begin{prop}\label{PFH1}
 Let ${\bf
 \bar{G}}$ and ${\bf G}$ be projectively related with the projective factor $P$.
  \ben
  \item[{\rm (i)}] (\cite{Shen6}) $\bar{R}^i_{\ k}=R^i_{\ k}$ if and only if $P$ is a Funk function
 of {\bf G}.
 \item[{\rm (ii)}] (\cite{Shen6}) $\bar{R}ic=Ric$ if and only if $P$ is a weak Funk function
 of {\bf G}.
 \item[{\rm (iii)}] If {\bf G} is of isotropic curvature, then ${\bf \bar{G}}$ is of isotropic
 curvature if and only if $P$ is a Hamel function of {\bf G}.
   \een
\end{prop}

{\it Proof :} Cases (i) and (ii)  directly follow from (\ref{pr9})
and (\ref{pr13}) respectively. For case (iii), it follows from
(\ref{pr9}) that the Riemann curvature tensor $\bar{R}^i_{\ k}$ of
${\bf \bar{G}}$ is given by
 \beq
 &&\hspace{2cm}\bar{R}^i_{\
 k}=\bar{R}\delta^i_k-\bar{\tau}_ky^i,\label{cwy76}\\
 &&\big(\bar{R}:=R+P^2-P_{;0},\ \ \ \bar{\tau}_k:=\tau_k+PP_k+P_{k;0}-2P_{;k}\big),\nonumber
 \eeq
 where we have put
 $$\bar{G}^i=G^i+Py^i,\ \ \ \ R^i_{\
k}=R\delta^i_k-\tau_ky^i.
 $$
Now by $R_{.i}=2\tau_i$, it  can be directly  verified that
 $$
\bar{R}_{.k}=2\bar{\tau}_k \Longleftrightarrow P_{.k;0}=P_{;k}.
 $$
 This completes the proof of Case (iii).  \qed

\begin{prop}\label{prop814}
Let ${\bf G}$ be a spray of isotropic curvature, and ${\bf
\bar{G}}$  be projectively related to {\bf G} with the projective
factor $P$ defined by  $P(y):=Q(-y)$ for a Hamel function $Q$ of
{\bf G}. Then ${\bf \bar{G}}$ is of isotropic
 curvature if and only if (\ref{cwy75}) holds.
\end{prop}

{\it Proof :} It directly follows from Proposition \ref{PFH1}
(iii) and Lemma \ref{Plem26}. \qed

\

A Funk metric on a convex domain of $R^n$ has a lot of special
curvature properties, and a Klein metric is produce by a Funk
metric. Using the Finsler spray induced by a Funk metric, we
define a family of sprays given by (\ref{ME1}), and we obtain
Theorem \ref{MTh1} on some properties of the spray (\ref{ME1}).
Now we give its proof.

\

{\it Proof of Theorem \ref{MTh1}:} Since
 $$
 \bar{G}^i=G^i+Py^i,\ \ \  G^i=\frac{1}{2}Fy^i,\ \ \ \
 \big(P(y):=cF(-y)\big),
 $$
 by ({\ref{cwy77}) and $F_{;i}=0$, we have
 \beq
 P_{;i}(y)&&\hspace{-0.6cm}=c\big[G^r_i(y)+G^r_i(-y)\big]F_{.r}(-y)\nonumber\\
 &&\hspace{-0.6cm}=-\frac{1}{2}cF(-y)F_{.i}(y)+c\big[\frac{1}{2}F(y)+F(-y)\big]F_{.i}(-y).\label{cwy79}
 \eeq

(i) Since {\bf G} is of isotropic (actually constant) curvature
and $F$ is a Hamel function of {\bf G}, it follows from
Proposition \ref{prop814} that ${\bf \bar{G}}$ is of isotropic
curvature with $\bar{R}^i_{\
k}=\bar{R}\delta^i_k-\frac{1}{2}\bar{R}_{.k}y^i$, where, by
(\ref{cwy76}) and ({\ref{cwy79}), $\bar{R}$ is written as
 \be\label{cwy78}
 \bar{R}(y)=-\frac{1}{4}F^2(y)+c^2F^2(-y)+c\big[F(y)F(-y)+F^2(-y)\big].
 \ee
By (\ref{cwy78}) and ({\ref{cwy79}), a direct computation gives
 \beq\label{cwy80}
 &&\hspace{1cm}\bar{R}_{\bar{;}i}=\bar{R}_{;i}-2\bar{R}P_i-P\bar{R}_{.i}=c(1+2c)F(-y)T_i,\\
 &&\Big(T_i:=2(1+c)F(-y)F_{.i}(-y)+2F(y)F_{.i}(-y)-F(-y)F_{.i}(y)\Big).\nonumber
 \eeq
If $\bar{R}=0$ or
$$
 y^iT_i=-F(-y)\big[(2(1+c)F(-y)+3F(y)\big]=0,
 $$
 we have $F(-y)=\sigma F(y)$ for some constant $\sigma$. By the
  geometric meaning, $F(-y)=\sigma F(y)$
 holds only on a nowhere dense set $\mathcal{A}$ of  $T\Omega\setminus
 \{0\}$. If ${\bf \bar{G}}$ is of constant curvature, we
 have $\bar{R}_{\bar{;}i}=0$ on the open set $\widetilde{\mathcal{A}}:=T\Omega\setminus
 \{0\}-\mathcal{A}$ since $\bar{R}\ne 0$ on
 $\widetilde{\mathcal{A}}$. There is some $T_i$ which is not  zero on an
 $\widetilde{\mathcal{A}}$. Then by (\ref{cwy80}), we have $c=0$
or $c=-1/2$. Conversely, if $c=0$ or $c=-1/2$, then
$\bar{R}_{\bar{;}i}=0$ by (\ref{cwy80}). So ${\bf \bar{G}}$ is of
constant curvature.

 (ii) Using the facts that ${\bf\bar{G}}$ is projectively flat and
 ${\bf\bar{G}}$ is of isotropic curvature by (i),
if ${\bf\bar{G}}$ is metrizable (on a conical region
$\mathcal{C}(\Omega)$),  we have $\bar{R}_{\bar{;}i}=0$ (on
$\mathcal{C}(\Omega)$) according to a result in \cite{Yang2}. Then
the above proof shows $c=0$ or $c=-1/2$.

Conversely,  if $c=0$, {\bf G} is induced by $F$. If $c=-1/2$,
then by (\ref{cwy78}), we obtain
 $$
 -\bar{R}=\frac{1}{4}\big[F(y)+F(-y)\big]^2,
 $$
which is a Finsler metric (called the Klein metric). Further by
(\ref{cwy80}) we have $\bar{R}_{\bar{;}i}=0$. Therefore, {\bf G}
is induced by the Klein metric $-\bar{R}$ according to a result in
\cite{Yang2}.     \qed

\section{Constructions of Hamel-Funk Functions}

 In this section, we consider the existence of Hamel functions
 or Funk functions on some special spray manifolds.

For a Finsler spray  {\bf G} induced by a Finsler metric $F$, it
is clear that $F$ is a Hamel function of {\bf G}, but it is not a
Funk function of {\bf G}. For a Minkowski spray ${\bf G}=0$, we
can construct many Funk functions {\bf G} by solving $P$ from
(\ref{prof1}) for a given $\varphi$. Then by Lemma \ref{HF2}, we
can obtain a lot of Hamel Function of ${\bf G}=0$. By Lemmas
\ref{Plem26}--\ref{Plem30}, we have different ways of constructing
Hamel functions of  a spray {\bf G}, especially when ${\bf G}=0$.
Theses Hamel functions include projectively flat Finsler metrics
of constant flag curvature.  We know that  a Hamel function is
projectively invariant (Lemma  \ref{PPHF1}); and all Hamel
functions of a spray form a vector space (Lemma \ref{Plem25}), but
it might be difficult to know its dimension.

A Hamel function or a Funk function might be also a Finsler
metric. A Hamel (resp. Funk) function is called a Finsler-Hamel
(resp. Finlser-Funk) function, if it is also a Finsler metric. The
following example gives a Finlser-Hamle function of a spray.

\begin{ex}\label{EEHF01}
 For an $(\alpha,\beta)$-metric $F=\alpha \phi(s),\  s=\beta/\alpha$,  the spray coefficients $\bar{G}^i$ of $F$
 is given by (cf. \cite{Shen7})
  \beqn
  &&\bar{G}^i=G^i+\alpha^{-1}\Theta (-2\alpha Q
  s_0+r_{00})y^i,\\
  &&G^i:=G^i_{\alpha}+\alpha Q s^i_{\ 0}+\Psi (-2\alpha Q s_0+r_{00})b^i,
  \eeqn
where $G^i_{\alpha}$ denote the spray coefficents of $\alpha$ and
 \beqn
  &&Q:=\frac{\phi'}{\phi-s\phi'},\ \
  \Theta:=\frac{Q-sQ'}{2\Delta},\ \
  \Psi:=\frac{Q'}{2\Delta},\ \ \Delta:=1+sQ+(b^2-s^2)Q',\\
  &&r_{ij}:=\frac{1}{2}(b_{i|j}+b_{j|i}),\ \ s_{ij}:=\frac{1}{2}(b_{i|j}-b_{j|i}),\ \
  s_j:=b^is_{ij},\ \ s^i_{\ j}:=a^{im}s_{mj},
 \eeqn
where  $b_{i|j}$ denotes the covariant derivative of the
Levi-Civita connection of $\alpha$.

It is clear that $F$ is a Finsler-Hamel function of the spray
$G^i$, since $F$ is naturally a Hamel function of $\bar{G}^i$ and
a Hamel function is projectively invariant (Lemma  \ref{PPHF1}).
\end{ex}

\begin{prop}\label{PEHF001}
 A  1-form is a Hamel function of a
 spray if and only if it is closed. So for a spray, a Hamel function plus a closed 1-form is a
 Hamel function.
\end{prop}

\begin{prop}\label{PEHF002}
 If {\bf G} is of isotropic curvature, then the
 $S$-curvature of {\bf G} with respect to any volume form is a Hamel function of {\bf G}.
\end{prop}

Proposition \ref{PEHF001} follows directly from the definition of
a Hamel function. For Proposition \ref{PEHF002}, the S-curvature
satisfies $S_{;i}-S_{.i;0}=0$ if {\bf G} is of isotropic curvature
(\cite{Shen2}). So the
 $S$-curvature of {\bf G} with respect to any volume form is a Hamel function of {\bf G}.

\

Now we turn to study Funk functions of certain sprays. We have
shown the construction of Funk functions of a Minkowski spray
${\bf G}=0$. Actually, we can prove that there exist local Funk
functions on a R-flat spray manifold, as shown in Theorem
\ref{MTh2}.

\

{\it Proof of Theorem \ref{MTh2} :} Consider the following PDE
system
 \be\label{EHF01}
Q_{;i}=QQ_{.i},\ \ \ \ \ \ YP=P,\ \ \ \  (Y:=y^i\dot{\pa}_i).
 \ee
We only need to prove that (\ref{EHF01}) is integrable if {\bf G}
is R-flat. The integrable conditions of (\ref{EHF01}) are
 \be\label{pr49}
 \delta_iYP-Y\delta_iP=[\delta_i,Y]P,\hspace{1.5cm}
 \delta_i\delta_jP-\delta_j\delta_iP=[\delta_i,\delta_j]P.
 \ee
The first condition in (\ref{pr49}) automatically holds, since by
(\ref{EHF01}), we have
$$[\delta_i,Y]=0,\ \ \delta_iYP=P_{;i},\ \ \
Y\delta_iP=Y\big(PP_{.i}\big)=P_{;i}.$$
 The second condition in
(\ref{pr49}) is equivalent to
 \be\label{EHF02}
 Q_{;i;j}-Q_{;j;i}=-Q_{.r}R^r_{\ ij}.
 \ee
Using (\ref{EHF01}), we have
 $$
 Q_{;i;j}=Q_{;j}Q_{.i}+QQ_{.i;j}.
 $$
Then using $Q_{.i;j}=Q_{;j.i}$ and (\ref{EHF01}), we get
 $$
Q_{;i;j}-Q_{;j;i}=0.
 $$
 Thus (\ref{EHF02}) holds if {\bf G} is R-flat ($R^r_{\ ij}=0$).
 \qed

 \

A similar system of PDEs to (\ref{EHF01}) is concerned in
\cite{Cram}. Next we consider the construction of Funk functions
of sprays with non-zero Riemann curvature, and prove some results
on the non-existence of non-zero Funk function on some spray
manifolds under certain conditions.

\begin{thm}(\cite{Shen6})\label{TEHF1}
A weak Funk function  on a complete regular spray manifold must be
zero.
\end{thm}

{\it Proof :} Let $({\bf G},M)$ be a complete regular spray
manifold, and $Q$ be a weak Funk function of $({\bf G},M)$. Then
$Q$ satisfies $Q_{;0}=Q^2$. Along a geodesic $c=c(t)$ of {\bf G}
with $c(0)=x\in M,\dot{c}(0)=y \in T_xM$, we have
 $$
 Q'(t)-Q^2(t)=0,\ \ \ \ (Q(t):=Q(c(t),\dot{c}(t))),
 $$
whose solution is given by
 $$
Q(t)=\frac{Q(0)}{1-Q(0)t}.
 $$
Since $({\bf G},M)$ is complete, $Q(t)$ should be defined on
$(-\infty,+\infty)$. It is easy to conclude $Q(0)=0$ since $Q(t)$
is unbounded at the point $t=1/Q(0)$ if $Q(0)\ne0$. Thus we get
$Q=0$  since the geodesic $c$ is arbitrarily given. \qed

\

In Theorem \ref{TEHF1}, if $({\bf G},M)$  is not complete, then
possibly there exist nontrivial Funk functions. For example,
consider the Minkowski spray ${\bf G}=0$ on a convex domain
$\Omega$ of $R^n$, and  the Funk metric on $\Omega$ is a Funk
function of {\bf G}.

A Funk functions on a spray manifold of scalar curvature has some
special properties. First we introduce a basic proposition as
 follows.

\begin{prop}\label{PEHF1}
Let {\bf G} be a spray of scalar curvature $R^i_{\
k}=R\delta^i_k-\tau_ky^i$ and $Q$ be a Funk function of {\bf G}.
Then $R,\tau_k$ and $Q$ are related by
 \be\label{EHF1}
 RQ_{.k}=Q\tau_k.
 \ee
  \ben
 \item[{\rm (i)}]  If {\bf G} is of isotropic curvature with $R\ne 0$, then (\ref{EHF1})
 gives
  \be\label{EHF2}
 Q^2=cR,
  \ee
  where $c=c(x)$ is a scalar function.
   \item[{\rm (ii)}]  If  {\bf G} is of constant curvature with $R\ne
  0$, then $Q$ is a closed 1-form. In this case, if $Q\ne 0$, then for some
  scalar function $\sigma=\sigma(x)$, we have
   \be\label{EHF3}
 Q=\sigma_0,\ \ \ \ R=\pm \big[(e^{-\sigma})_{0}\big]^2,\ \ \ \
 Q^2=\pm e^{2\sigma}R.
   \ee
  \een
\end{prop}

{\it Proof :} Since $Q$ is a Funk function of {\bf G}, we have
$Q_{.r}R^r_{\ k}=0$ by (\ref{EHF02}). Using $R^i_{\
k}=R\delta^i_k-\tau_ky^i$ we get (\ref{EHF1}) immediately.

If {\bf G} is of isotropic curvature with $R\ne 0$, then plugging
$R_{.k}=2\tau_k$ into (\ref{EHF1}) yields
 $$
 \frac{2Q_{.i}}{Q}=\frac{R_{.i}}{R}, \ \ \ (letting \ Q\ne 0),
 $$
which gives (\ref{EHF2}) by an integration.

If {\bf G} is of constant curvature with $R\ne 0$, then we have
$R_{;i}=0$. Thus by (\ref{EHF2}), we have
 $$
 2Q^2Q_{.i}=2QQ_{;i}=Rc_i+cR_{;i}=Rc_i,\ \ \ (c_i:=c_{x^i}).
 $$
 Contracting the above by $y^i$ and using (\ref{EHF2}), we have
  $$
2cQ=c_0.
  $$
If $c=0$, then $Q=0$ by (\ref{EHF2}). If $c\ne0$, then we have
$$Q=(\ln\sqrt{|c|})_0.$$
 So $Q$ is a closed 1-form. Putting $Q=\sigma_0$ with $c=\pm e^{2\sigma}$, we obtain (\ref{EHF3}). \qed

\begin{thm} (\cite{BM0})\label{TEHF2}
 A Funk  function on  a Finsler manifold of non-zero scalar flag curvature
 must be zero. So a spray {\bf G} of scalar curvature with non-zero
 Riemann curvature is not metrizable if there is a non-zero Funk
 function of {\bf G}.
\end{thm}

{\it Proof :} Let $L$ be a Finsler metric of scalar flag curvature
$\lambda\ne 0$, and $Q$ be a Funk metric of {\bf G}. Then the
spray {\bf G} of $L$ is given by
 $$
R^i_{\ k}=R\delta^i_k-\tau_ky^i,\ \ \ (R:=\lambda L,\ \
\tau_k:=\lambda y_k).
 $$
Suppose $Q\ne 0$. Then by (\ref{EHF1}) we have
 $$
 \frac{2Q_{.i}}{Q}=\frac{L_{.i}}{L},
 $$
which gives $Q^2=cL$ for a scalar function $c=c(x)\ne 0$. By
$Q^2=cL$, we have
 $$
 2Q^2Q_{.i}=2QQ_{;i}=Lc_i,\ \ \ (c_i:=c_{x^i}).
 $$
 Contracting the above by $y^i$ and using $Q^2=cL$, we see that
 $Q=c_0/(2c)$ is a closed 1-form. Thus  $L=Q^2/c$ is not a Finsler metric, which gives a contradiction.  \qed

\

In the following examples, we consider the construction of Funk
functions of  a projectively flat Berwald spray of isotropic
curvature. It should be the first time that we obtain non-trivial
Funk functions on a spray manifold with non-zero Riemann
curvature.

\begin{ex}\label{EEHF1}
 Let {\bf G} be a projectively flat Berwald spray with
  \be\label{EHF4}
 G^i:=\tau_0y^i,\ \ \ \ \big(\tau_i:=\tau_{x^i},\ \ \tau:=f(\sigma)\big),
  \ee
  where $\sigma=\sigma(x)$ is a scalar function.
  Suppose
  $Q:=\sigma_0$  is a Funk function of
  {\bf G}. Then by $Q_{;j}=QQ_{.j}$, we have
   $$
 \sigma_{0;j}-\sigma_0\sigma_j=0.
   $$
  which can be rewritten  as
    $$
 \sigma_{ij}-\big[1+2f'(\sigma)\big]\sigma_{i}\sigma_{j}=0,\ \ \ \
 or\ \ \
  \big[\int e^{-\sigma-2f(\sigma)}d\sigma\big]_{ij}=0.
    $$
    Thus $\sigma$ is given by
     \be\label{EHF5}
 \int e^{-\sigma-2f(\sigma)}d\sigma=\langle a,x\rangle+b,
     \ee
     where $a$ is an $n$-vector and $b$ is a constant.

     (i) By the different choices of the function $f$ in (\ref{EHF5}), we
     obtain many Funk functions $Q=\sigma_0$ of the spray {\bf G}
     in (\ref{EHF4}).

     (ii) {\bf G} in (\ref{EHF4}) is of
     isotropic curvature $R^i_{\ k}=R\delta^i_k-\frac{1}{2}R_{.k}y^i$
  with
   \be\label{EHF04}
   R:=-(\tau_0)^2-\tau_{0;0}=-\big\{f''(\sigma)+f'(\sigma)[1+f'(\sigma)]\big\}(\sigma_0)^2.
  \ee
  It follows from Proposition \ref{PEHF1}(i)
     that (\ref{EHF2}) holds, and then by (\ref{EHF04}),  $c$ in (\ref{EHF2}) is given by
      \be\label{EHF6}
 c=-\big\{f''(\sigma)+f'(\sigma)[1+f'(\sigma)]\big\}^{-1}
      \ee

(iii) By Theorem \ref{TEHF2}, the spray in (\ref{EHF4}) with
$\sigma$ given by (\ref{EHF5}) is non-metrizable.
\end{ex}

In Example \ref{EEHF1}, taking $f(t):=t$ as a special case, by
(\ref{EHF4}), (\ref{EHF5}) and (\ref{EHF6}),
      we obtain
       $$
 G^i=-\frac{1}{3}\frac{\langle a,y\rangle}{\langle a,x\rangle+b}y^i,\ \ \ \
 Q=-\frac{1}{3}\frac{\langle a,y\rangle}{\langle a,x\rangle+b},\ \ \ \ c=-\frac{1}{2}.
       $$

\begin{ex}\label{EEHF2}
 If the spray {\bf G} given by (\ref{EHF4}) is of constant
 curvature, we can also obtain nontrivial Funk function of {\bf
 G}.  Let {\bf G} be given by (\ref{EHF4}) and $\sigma$ is
 determined by (\ref{EHF5}). Then $Q=\sigma_0$ is a Funk function
 of {\bf G} by Example \ref{EEHF1}. Now we can choose $f$ such that  {\bf G} is of constant curvature.
 To do so, we need to solve $R_{;i}=0$ for  $R$ given by
 (\ref{EHF04}). It is easy to see that $R_{;i}=0$ if and  only if
  $$
 f'''(t)+3f''(t)+2f'(t)f''(t)+2f'(t)+2(f'(t))^2=0,
  $$
  whose solution is given by
   $$
 f(t)=-c_1e^{-t}+\ln\big|1+c_2e^{2c_1e^{-t}}\big|+c_3,
   $$
   where $c_1,c_2,c_3$ are constant.
\end{ex}

\begin{ex}\label{EEHF3}
 Consider an $n$-dimensional projectively flat Berwald spray {\bf G}:
  $$
 G^i=-\frac{\langle x,y\rangle}{|x|^2}y^i.
  $$
 In \cite{Yang2}, we have shown that {\bf G} is of constant curvature and it is not
 metrizable, and
 further, the Riemann curvature of {\bf G} satisfies
  \be\label{EHF7}
 R^r_{\ jk}=\frac{1}{2}R_{.j}\delta^r_k-(jk),\ \ \ \ R:=\frac{|x|^2|y|^2-\langle
 x,y\rangle^2}{|x|^4}.
  \ee
   We will prove that there is a Funk function
  $Q\ne 0$ of {\bf G} if and only if $n=2$ and $Q$ is given by
   \be\label{EHF8}
 Q=\sigma_0,\ \  \ \
 \sigma:=\ln\Big|\frac{c_1}{c_1\arctan\frac{x^2}{x^1}+c_2}\Big|,
   \ee
   where $c_1,c_2$ are constant. In this case, $Q,R$ satisfy (\ref{EHF3}), and  $c$ in (\ref{EHF2}) is given by
      \be\label{EHF08}
  c=\Big(\frac{c_1}{c_1\arctan\frac{x^2}{x^1}+c_2}\Big)^2.
      \ee

       Actually, by Proposition
   \ref{PEHF1}(ii), we first have $Q=\sigma_0$ for some scalar function
   $\sigma=\sigma(x)$. By (\ref{EHF02}), the integral condition
   for $Q$ is given by $Q_{.r}R^r_{\ jk}=0$. Then by (\ref{EHF7}),
   the integral condition becomes
    \be\label{EHF9}
 (|x|^2\delta_{ij}-x^ix^j)\sigma_k-(jk)=0.
    \ee
    Taking a sum over $i,j$ in (\ref{EHF9}) we have
    \be\label{EHF10}
 (n-2)|x|^2\sigma_k+\sigma_rx^rx^k=0.
    \ee
    Then a further contraction in (\ref{EHF10}) by $x^k$ gives
    $\sigma_rx^r=0$, which implies that $\sigma$ is constant by (\ref{EHF10})
    again when $n>2$. When  $n=2$,  the integral condition (\ref{EHF9}) becomes $x^1\sigma_1+x^2\sigma_2=0$,
    and in this case, solving the following system of PDEs
     $$
  Q_{;i}=QQ_{.i}\ \ \ \ \ \Big(\Longleftrightarrow
  \sigma_{ij}+\frac{1}{|x|^2}(\sigma_ix^j+\sigma_jx^i)-\sigma_i\sigma_j=0\Big),
     $$
     we obtain  (\ref{EHF8}). By the
     way, we also have $n=2$ by the $R$ in (\ref{EHF7}) and
     (\ref{EHF3}).
     .
\end{ex}

\section{ Hamel-Funk Sprays}

In this section  we introduce a new spray associated to a spray
and a Hamel function, and then study some properties of this new
spray.

For a given spray {\bf G} and a nonzero homogeneous function $Q$
of degree one, we consider a spray ${\bf \bar{G}}$ defined by
 \be\label{HFS1}
\bar{G}^i:=G^i+\frac{Q_{;0}}{2Q}y^i.
 \ee
When ${\bf G},Q$ satisfies certain good conditions, ${\bf
\bar{G}}$ is metrizable (Theorem \ref{THF1} below). So in general
case, ${\bf \bar{G}}$ is pre-metrizable. For short, ${\bf
\bar{G}}$ is called a PM-spray, and ${\bf G},Q$ are respectively
called the adjoint spray and adjoint function of ${\bf \bar{G}}$.

\begin{lem}\label{LHFS001}
Each spray in the projective class $Proj({\bf G})$ of a spray {\bf
G} is a local PM-spray ${\bf \bar{G}}$ defined by (\ref{HFS1}) for
some $Q$.
\end{lem}

{\it Proof :} For an arbitrary spray ${\bf \bar{G}}\in Proj({\bf
G})$, we have $\bar{G}^i=G^i+Py^i$ for some function $P$. Consider
the following system of PDEs with the unknown function $Q$:
 \be\label{HFS001}
 \bar{Y}(Q)=2PQ,\ \ \ \ \ Y(Q)=Q,\ \ \ \  \big(\bar{Y}:=y^i\delta_i,\ \
 Y=y^i\dot{\pa} y^i\big).
 \ee
The integral condition of (\ref{HFS001}) is given by
 $$
 [\bar{Y},Y](Q)=\bar{Y}Y(Q)-Y\bar{Y}(Q).
 $$
Applying $[\bar{Y},Y]=-\bar{Y}$ and (\ref{HFS001}), it is easy to
verify that the above integral condition is satisfied. So
(\ref{HFS001}) has a local solution $Q$.    \qed

\

In Lemma \ref{LHFS001}, $Q$ is generally not unique. For instance,
letting
 \be\label{HFS002}
 G^i:=0,\ \ \ \ \ P:=-\frac{x^1y^1+x^2y^2}{1+(x^1)^2+(x^2)^2},\ \ \ \ \
 Q:=\frac{f(\frac{y^2}{y^1},\frac{x^2y^1-x^1y^2}{y^1})}{1+(x^1)^2+(x^2)^2}\ y^1,
 \ee
 where $f=f(s,t)$ is some function.
Then $Q$ satisfies (\ref{HFS001}).

\begin{lem}\label{LHFS01}
 A  PM-spray  is projectively
 invariant with  a fixed adjoint  function.
\end{lem}

{\it Proof :} Let ${\bf \bar{G}}_1, {\bf \bar{G}}_2$ be two
PM-sprays associated to $({\bf G}_1,Q),({\bf G}_2,Q)$
respectively, where ${\bf G}_1$ and ${\bf G}_2$ are projectively
related with $G_1^i=G_2^i+Py^i$. By a direct computation, we have
(see (\ref{proj1}))
 $$
 Q_{|i}=Q_{;i}-(PQ)_{.i},\ \ \ \ Q_{|0}=Q_{;0}-2PQ,
 $$
where $_|$ and $_{;}$ denote the horizontal derivatives of ${\bf
G}_1$ and ${\bf G}_2$ respectively. Therefore, we obtain
 \beqn
 \bar{G}^i_1&&\hspace{-0.6cm}=G^i_1+\frac{Q_{|0}}{2Q}y^i=G_2^i+Py^i+\frac{Q_{;0}-2PQ}{2Q}y^i\\
&&\hspace{-0.6cm}=G^i_2+\frac{Q_{;0}}{2Q}y^i=\bar{G}^i_2,
 \eeqn
which completes the proof. \qed

\

Lemma \ref{LHFS01} is similar to the result of the spray defined
by the S-curvature (\cite{Shen2}). For a PM-spray associated to
$({\bf G},Q)$, one important case is when $Q$ is a Hamel or Funk
function of {\bf G}. We
 have the following definition.

\begin{Def}\label{DHFS63}
Let {\bf G} be a spray and $Q$ be a Hamel function of {\bf G}. The
spray ${\bf \bar{G}}$ defined by  (\ref{HFS1}) is called a Hamel
spray associated to $({\bf G},Q)$, and ${\bf G},Q$ are
respectively called the adjoint spray and adjoint Hamel function
of ${\bf \bar{G}}$.
 As a special case, if $Q$ is a Funk function of {\bf
G}, then ${\bf \bar{G}}$ is called a Funk spray associated to
$({\bf G},Q)$.
\end{Def}

It is clear that the Funk spray ${\bf\bar{G}}$ associated to
$({\bf G},Q)$ has the form $\bar{G}^i=G^i+(Q/2)y^i$.

\begin{lem}\label{LHFS02}
A  Hamel spray  is projectively
 invariant with  a fixed adjoint Hamel function.
\end{lem}

{\it Proof :} It directly follows from Lemma \ref{LHFS01} and the
fact that a Hamel function is projectively invariant, and so $Q$
is a Hamel function of ${\bf G}_1$ if and only if it is  a Hamel
function of ${\bf G}_2$, where ${\bf G}_1$ and ${\bf G}_2$ are
projectively related.   \qed

\

 It is clear that a Finsler spray ${\bf \bar{G}}$ induced by a Finsler metric $\bar{F}$ is  a Hamel spray
 associated to $({\bf \bar{G}},\bar{F})$, and then by Lemma \ref{LHFS02}, ${\bf \bar{G}}$ is a Hamel spray
 associated to $({\bf G},\bar{F})$ for any spray {\bf G}
 projectively related to ${\bf \bar{G}}$.

\

{\it Proof of Theorem \ref{THF1} :}  If a spray {\bf G} is
projectively metrizable, then there is a Finsler metric $\bar{F}$
such that $\bar{G}^i:=G^i+Py^i$ for some $P$ is the Finsler spray
of $\bar{F}$. By Lemma \ref{Prelem2}(ii), $\bar{F}$ is a Hamel
function of {\bf G}. This fact can also be proved by  that
$\bar{F}$ is a natural Hamel function of ${\bf \bar{G}}$ and a
Hamel function is projectively invariant. In this case, $P$ is
given by (\ref{pr4}), and so the spray ${\bf \bar{G}}$ of
$\bar{F}$ is a Hamel spray associated to $({\bf G},\bar{F})$.

Conversely, if a spray {\bf G}  has a Finsler-Hamel function
$\bar{F}$, it follows from Lemma \ref{Prelem2} that {\bf G} is
projectively  related to the Finsler spray ${\bf\bar{G}}$ of
$\bar{F}$ which is the Hamel spray associated to $({\bf
G},\bar{F})$. This finishes the proof.           \qed

\

By Theorem \ref{THF1} we have the following corollary.

\begin{cor}
 If a spray {\bf G} has a Finsler-Funk function $\bar{F}$, then {\bf G}  is projectively induced
by  $\bar{F}$ whose spray is the Funk spray
 associated to $({\bf G},\bar{F})$.
\end{cor}

\begin{ex}
 In Example \ref{EEHF01}, $\bar{F}$ is a Finsler-Hamel function of
 ${\bf \bar{G}}$. So by Theorem \ref{THF1},  ${\bf \bar{G}}$
 is a Hamel spray associated to $({\bf G},\bar{F})$.
\end{ex}

For a Hamel or Funk spray, we give some special properties as
follows.

\begin{lem}\label{LHFS1}
Let ${\bf \bar{G}}$ and  {\bf G} be two sprays  projectively
related and $Q\ne 0$ be a Hamel function. Then  ${\bf \bar{G}}$ is
a Hamel spray associated to $({\bf G},Q)$ if and only if
  $Q_{\bar{;}k}=0$.
\end{lem}

{\it Proof :} Put $\bar{G}^i=G^i+Py^i$. A direct computation gives
 \be\label{HFS37}
Q_{\bar{;}i}=Q_{;i}-Q_{.r}(P_{.i}y^r+P\delta^r_i)=Q_{;i}-(PQ)_{.i}.
 \ee

If ${\bf \bar{G}}$ is a Hamel spray associated to $({\bf G},Q)$,
then by definition we have
 \be\label{HFS38}
  P=\frac{Q_{;0}}{2Q},
  \ee
where $Q$ is a Hamel function of the spray {\bf G}.  It follows
from (\ref{HFS37}) that
 \beqn
Q_{\bar{;}i}=Q_{;i}-\frac{1}{2}Q_{;0.i}=\frac{1}{2}(Q_{;i}-Q_{.i;0})=0.
 \eeqn

Conversely, if $Q_{\bar{;}k}=0$,  then since $Q$ is a Hamel
function, it follows from (\ref{HFS37}) that
 $$
 0=Q_{\bar{;}i}=Q_{.i;0}-(PQ)_{.i}=\frac{1}{2}Q_{;0.i}-(PQ)_{.i}=(\frac{1}{2}Q_{;0}-PQ)_{.i},
 $$
 which implies (\ref{HFS38}). Thus ${\bf \bar{G}}$ is
a Hamel spray associated to $({\bf G},Q)$.
 \qed

\begin{thm}\label{THFS2}
(i) An $n$-dimensional Hamel spray ${\bf \bar{G}}$    associated
to $({\bf G},Q)$ is of scalar curvature if and only if its
Riemmann curvature has the form
 \be\label{HFS4}
 \bar{R}^i_{\ k} =\lambda (Q^2\delta^i_k-QQ_{.k}y^i),
 \ee
 where $\lambda=\lambda(x,y)$ is a scalar function given by
  \be\label{HFS01}
\lambda:=(R+P^2-P_{;0})/Q^2,\ \ \ \ \big(P:=Q_{;0}/(2Q),\ \
R:=Ric/(n-1)\big).
  \ee

 (ii) A Funk spray ${\bf \bar{G}}$    associated to $({\bf G},Q)$ is of
isotropic curvature if and only if its Riemmann curvature has the
from (\ref{HFS4})  with $\lambda=\lambda(x)$ being given by
 \be\label{HFS02}
 \lambda=\eta-\frac{1}{4},\ \ \ \big(\eta:=R/Q^2\big).
 \ee
\end{thm}

{\it Proof :}  Put
 $$
\bar{G}^i=G^i+Py^i,\ \ \  \big(P:=Q_{;0}/(2Q)\big).
 $$
Let  ${\bf \bar{G}}$ be of scalar curvature. Then {\bf G} is also
of scalar curvature and we put $R^i_{\ k}=R\delta^i_k-\tau_ky^i$.
Then we have (see (\ref{cwy76}))
 \beq
 &&\hspace{2cm}\bar{R}^i_{\
 k}=\bar{R}\delta^i_k-\bar{\tau}_ky^i,\label{HFS3}\\
 &&\big(\bar{R}:=R+P^2-P_{;0},\ \ \ \bar{\tau}_k:=\tau_k+PP_{.k}+P_{.k;0}-2P_{;k}\big).\nonumber
 \eeq
It is easy to obtain
 \beqn
&&P_{;k}=\frac{QQ_{;0;k}-Q_{;0}Q_{;k}}{2Q^2},\ \ \ \ \ \
P_{;0}=\frac{QQ_{;0;0}-(Q_{;0})^2}{2Q^2},\\
&&P_{.k}=\frac{2QQ_{;k}-Q_{;0}Q_{.k}}{2Q^2},\ \ \ \ \ \
P_{.k;0}=\frac{Q_{;k;0}}{Q}-\frac{3Q_{;0}Q_{;k}+Q_{;0;0}Q_{.k}}{2Q^2}+\frac{(Q_{;0})^2Q_{.k}}{Q^3},
 \eeqn
where in the third formula,  $Q_{;0.k}=2Q_{;k}$ is used, and in
the fourth formula, $Q_{.k;0}=Q_{;k}$ is concerned. Plugging   the
above formulas into (\ref{HFS3})  and applying
 $$
 Q_{;0;k}-Q_{;k;0}=-Q_{.r}R^r_{\ k},\ \ \ \ \ R^r_{\
 k}=R\delta^r_k-\tau_ky^r,
 $$
we obtain
 $$
 \bar{R}=\lambda Q^2,\ \ \ \ \bar{\tau}_k=\lambda QQ_{.k},
 $$
where $\lambda$ is given by (\ref{HFS01}). This gives the proof of
case (i).

If ${\bf \bar{G}}$ is of isotropic curvature and $Q$ is a Funk
 function of {\bf G}, then by Proposition \ref{PFH1}(iii), we see that ${\bf
 G}$ is of isotropic curvature, since $P=Q/2$ is a Hamel
 function. By (\ref{P315}) we get
  $$
 \frac{P^2-P_{;0}}{Q^2}=-\frac{1}{4}.
  $$
By (\ref{EHF2}), $\eta:=R/Q^2$ is just a function of $x$. So
(\ref{HFS02}) follows from (\ref{HFS01}).   \qed

\

Based on Theorem \ref{THFS2}, we make the following definition.

\begin{Def}
 A Hamel spray ${\bf \bar{G}}$    associated to $({\bf
 G},Q)$ is called of scalar curvature $\lambda$ if its Riemann curvature has the form
 (\ref{HFS4}). ${\bf \bar{G}}$ is called of isotropic (resp.
 constant) curvature $\lambda$ if $\lambda=\lambda(x)$ (resp.
 $\lambda=constant$).
\end{Def}

Write (\ref{HFS4}) in the form
 $$
\bar{R}^i_{\ k} =\bar{R}\delta^i_k-\bar{\tau}_ky^i,\ \ \ \
\bar{R}:=\lambda Q^2,\ \ \ \  \bar{\tau}_k:=\lambda QQ_{.k}.
 $$
Then  it is easily to show that for $Q\ne 0$,
 $$
\lambda=\lambda(x)\ \Longleftrightarrow \
\bar{R}_{.k}=2\bar{\tau}_k;\ \ \ \ \lambda=constant\
\Longleftrightarrow \ \bar{R}_{\bar{;}k}=0,
 $$
where we have applied Lemma \ref{LHFS1} and Lemma \ref{Yplem} for
the second case in the above.

\begin{ex}\label{EHFS2}
Define a spray {\bf G} and 1-form $Q$ by
 $$
 G^i:=f'(\sigma)\sigma_0y^i,\ \ \ \ Q:=\sigma_0,
 $$
 where $f$ is a function on $R^1$ and $\sigma=\sigma(x)$ is given
 by (\ref{EHF5}). Then by Example \ref{EEHF1}, {\bf G} is of
 isotropic curvature and $Q$ is a Funk function of {\bf G}. By
 Theorem \ref{THFS2}(ii), the Funk spray ${\bf \bar{G}}$
 associated to $({\bf G},Q)$ satisfies
  \beqn
 \bar{R}^i_{\ k} =\lambda
 (Q^2\delta^i_k-QQ_{.k}y^i),\ \ \ \ \lambda:=-f''(\sigma)-f'(\sigma)[1+f'(\sigma)]-\frac{1}{4},
  \eeqn
  where we have have applied (\ref{HFS02}) and (\ref{EHF6}). So ${\bf \bar{G}}$ is of isotropic curvature $\lambda$. For
  suitable choice of $f$, $\lambda$ in the above can be a constant
  (cf. Example \ref{EEHF1}). For instance, taking $f(t)=t$ we get
  $\lambda=-9/4$.
\end{ex}

\begin{ex}\label{EHFS3}
Define a two-dimensional spray {\bf G} and 1-form $Q$ by
 $$
 G^i=-\frac{x^1y^1+x^2y^2}{(x^1)^2+(x^2)^2}y^i,\ \ \ \ \  Q:=\sigma_0,\ \  \ \
 \sigma:=-\ln\big|\kappa+\arctan\frac{x^2}{x^1}\big|,
  $$
  where $\kappa$ is a constant. By Example \ref{EEHF3}, {\bf G} is
  of constant curvature and $Q$ is a Funk function of {\bf G}. By
 Theorem \ref{THFS2}(ii), the Funk spray ${\bf \bar{G}}$
 associated to $({\bf G},Q)$ satisfies
  \beqn
 \bar{R}^i_{\ k} =\lambda
 (Q^2\delta^i_k-QQ_{.k}y^i),\ \ \ \ \lambda:=\big(\kappa+\arctan\frac{x^2}{x^1}\big)^2-\frac{1}{4},
  \eeqn
  where we have applied (\ref{HFS02}) and (\ref{EHF08}). So ${\bf \bar{G}}$ is of isotropic curvature $\lambda$.
\end{ex}

By Example \ref{EHFS3}, the  Funk spray ${\bf \bar{G}}$
 associated to $({\bf G},Q)$ is not necessarily of constant curvature even if {\bf G} is of constant
 curvature.

 \

For a spray {\bf G}, define a map
 \be\label{HHFS}
\mathcal{H}:[Q]\rightarrow \frac{Q_{;0}}{2Q},
 \ee
where $Q\ne0$ is a homogeneous function of degree one and $[Q]$
denotes  the equivalent class consisting of all the functions $cQ$
for arbitrary nonzero constant $c$.

Let $\mathcal{A}$ be the set of all nonzero Hamel functions and
set $\mathcal{H}_1:=\mathcal{H}|_{[\mathcal{A}]}$. By definition,
 we have $Q\in\mathcal{H}_1^{-1}(P)$ if and only if
   \be\label{HFS05}
 Q_{;i}=Q_{.i;0},\ \ \ \ \  Q_{;0}=2PQ.
 \ee
It is clear that for a given spray ${\bf G}$ and a function $P$,
the spray $\bar{G}^i:=G^i+Py^i$ is metrizable if and only if there
is a Finsler metric $Q$ satisfying (\ref{HFS05}).

Generally, $\mathcal{H}$ is not injective (see the example given
by (\ref{HFS002})). However, for certain function $P$, there will
have a unique element in $\mathcal{H}_1^{-1}(P)$. We first show
the following lemma.

\begin{lem}\label{LHFS3}
 Let $Q$ be a Hamel function of a spray {\bf G} and $\lambda=\lambda(x,y)$ be a homogeneous function of degree zero.
  Then $\lambda$ satisfies  $\lambda_{;0}=0$ and $\bar{Q}:=\lambda Q$   is
 a Hamel function of {\bf G} if and only if
 \be\label{HFS04}
 2Q\lambda_{;i}=Q_{;0}\lambda_{.i},\ \ \ \ or\ \ \lambda_{\bar{;}i}=0
 \ee
 where in the second case, the derivative is take with respect to
 the spray
 $$\bar{G}^i:=G^i+\frac{Q_{;0}}{2Q}y^i.$$
\end{lem}

{\it Proof :} Let $\bar{Q}=\lambda Q$ with $\lambda_{;0}=0$. We
have
 \beqn
&& \bar{Q}_{;i}=Q\lambda_{;i}+\lambda Q_{;i},\ \ \ \
\bar{Q}_{.i}=Q\lambda_{.i}+\lambda
 Q_{.i},\\
 &&\bar{Q}_{.i;0}=Q\lambda_{.i;0}+Q_{;0}\lambda_{.i}+\lambda
 Q_{.i;0}=Q\lambda_{.i;0}+Q_{;0}\lambda_{.i}+\lambda
 Q_{;i},
 \eeqn
where in the third formula above we have used the assumptions
$\lambda_{;0}=0$ and $Q_{.i;0}=Q_{;i}$. Therefore, $\bar{Q}$ is a
Hamel function if and only if
 $$
Q\lambda_{;i}+\lambda
Q_{;i}=Q\lambda_{.i;0}+Q_{;0}\lambda_{.i}+\lambda
 Q_{;i},
 $$
which reduces to the first formula in (\ref{HFS04}) since
$\lambda_{.i;0}=\lambda_{;0.i}-\lambda_{;i}=-\lambda_{;i}$.
Conversely, if the first formula in (\ref{HFS04}) holds,  it is
easy to get $\lambda_{;0}=0$ by contracting this formula by $y^i$.
Further, the above proof has shown that $\bar{Q}$ is Hamel
function.

Now under the spray ${\bf \bar{G}}$, a direct computation gives
 $$
\lambda_{\bar{;}i}=\lambda_{;i}-\frac{Q_{;0}}{2Q}\lambda_{.i},
 $$
which implies $\lambda_{\bar{;}i}=0\Longleftrightarrow
2Q\lambda_{;i}=Q_{;0}\lambda_{.i}$.
  \qed

\begin{prop}\label{PHFS3}
For a given spray {\bf G} and a function $P$, let $Q_o$ be a
solution of (\ref{HFS05}). Suppose that the spray
$\bar{G}^i:=G^i+Py^i$ is of nonzero scalar curvature. Then any
solution $Q$ of (\ref{HFS05}) is given by $Q=cQ_o$ for some
constant $c$, or $\mathcal{H}_1^{-1}(P)$ has a unique element.
\end{prop}

{\it Proof :}  Let $Q=cQ_o$ for a homogeneous function $c=c(x,y)$
of degree zero. Then we have
 $$
 \frac{Q_{o;0}}{2Q_o}=\frac{Q_{;0}}{2Q}\ (=P),
 $$
which gives $c_{;0}=0$.
 Since $Q,Q_o$ are Hamel functions of
{\bf G}, by Lemma \ref{LHFS3}, we have $c_{\bar{;}i}=0$.   By
assumption, ${\bf \bar{G}}$ is of nonzero scalar curvature. Then
$c$ is a constant due to Lemma \ref{Yplem}, which also shows that
$\mathcal{H}_1^{-1}(P)$ has a unique element. \qed

\begin{ex}
 In the system (\ref{HFS05}), put
  $$
 G^i:=0,\ \ \ \  P:=-\frac{\mu\langle x,y\rangle}{1+\mu|x|^2},\ \ \
 \ (\mu=constant \ne 0).
  $$
  A direct verification shows that
   $$
 Q_o:=\frac{\sqrt{(1+\mu|x|^2)|y|^2-\mu\langle
 x,y\rangle^2}}{1+\mu|x|^2}
   $$
   is a solution of (\ref{HFS05}). The spray $\bar{G}^i:=G^i+Py^i$
   is a Hamel spray associated to $({\bf G},Q_o)$ and it is of
   constant curvature $\mu$. By Proposition \ref{PHFS3}, any
solution $Q$ of (\ref{HFS05}) is given by $Q=cQ_o$ for some
constant $c$. Actually, ${\bf \bar{G}}$ is metrizable induced by
the metric $Q_o$. This can be an example of Theorem \ref{THFS02}.
\end{ex}

{\it Proof of Theorem \ref{THFS02} :} Let ${\bf \bar{G}}$ be a
Hamel spray associated to $({\bf G},Q)$ and ${\bf \bar{G}}$ be of
nonzero scalar curvature (cf. (\ref{HFS4}) with $\lambda\ne 0$).
If ${\bf \bar{G}}$ is metrizable, then there is a Finselr metric
$\bar{F}$ inducing ${\bf \bar{G}}$, and $\bar{F}$ is a Hamel
function of {\bf G}. By Lemma \ref{Prelem2}, ${\bf \bar{G}}$ can
be expressed as
 $$
\bar{G}^i=G^i+\frac{\bar{F}_{;0}}{2\bar{F}}y^i,
 $$
which implies
 $$\bar{F}_{;0}=2P\bar{F},\ \ \ \ \big(P:=Q_{;0}/(2Q)\big).$$
  Thus $\bar{F}$ is a
solution of (\ref{HFS05}). By assumption,  $Q$ is also a solution
of (\ref{HFS05}). Then by Proposition \ref{PHFS3}, $Q=c\bar{F}$
for some constant $c \ne 0$ is a Finsler metric. The converse is
clear by Theorem \ref{THF1}. \qed

\

 The Funk sprays in Examples \ref{EHFS2} and \ref{EHFS3} are not
 metrizable by Theorem \ref{THFS02}, since the scalar curvature is
 not zero, where in Example \ref{EHFS2}, $f(t)\ne -t/2+\ln|a+bt|$
 for any constant $a,b$.
Since any closed 1-form is a Hamel function of any spray, we have
the following corollary by Theorem \ref{THFS02}.

\begin{cor}\label{CHFS001}
 Let $\sigma=\sigma(x)\ne 0$ be a scalar function, and define a spray ${\bf
 \bar{G}}$ by
  \be\label{HFS005}
 \bar{G}^i=\frac{\sigma_{00}}{2\sigma_{0}}y^i,\ \ \ \ \
 (\sigma_i:=\sigma_{x^i},\ \ \sigma_{ij}:=\sigma_{x^ix^j}).
  \ee
  Then ${\bf \bar{G}}$ is metrizable if and only if ${\bf
  \bar{G}}$ is of zero scalar curvature, or equivalently,
   \be\label{HFS06}
 3(\sigma_{00})^2=2\sigma_0\sigma_{000},\ \ \ (\sigma_{ijk}:=\sigma_{x^ix^jx^k}).
   \ee
\end{cor}

{\it Proof :} ${\bf \bar{G}}$ is a Hamel spray associated to
$({\bf G},\sigma_0)$ with ${\bf G}=0$, and ${\bf \bar{G}}$ is of
scalar curvature. By (\ref{HFS01}), ${\bf
  \bar{G}}$ is of zero scalar curvature if and only if
  (\ref{HFS06}) holds. If ${\bf
  \bar{G}}$ is of zero scalar curvature, then ${\bf
  \bar{G}}$ is (locally) metrizable due to \cite{Yang2}. If ${\bf
  \bar{G}}$ is  of nonzero scalar curvature, then ${\bf
  \bar{G}}$ is not metrizable by Theorem \ref{THFS02}, since
  $\sigma_0$ is a not a Finsler metric.    \qed

\noindent Guojun Yang \\
Department of Mathematics \\
Sichuan University \\
Chengdu 610064, P. R. China \\
yangguojun@scu.edu.cn

\end{document}